\documentclass[12pt]{amsart}
\usepackage{amssymb}
\usepackage{amsmath}
\usepackage{mathabx}
\usepackage{amsthm}
\usepackage{amstext}
\usepackage{amsopn}
\usepackage{mathrsfs}
\usepackage{latexsym}
\DeclareMathAlphabet{\mathpzc}{OT1}{pzc}{m}{it}
\allowdisplaybreaks
\newtheorem{Definition}{Definition}[section]

\newtheorem{Theorem}{Theorem}[section]
\newtheorem{Corollary}{Corollary}[section]
\newtheorem{Remark}{Remark}[section]
\newtheorem{Example}{Example}[section]

\begin{document}
\bibliographystyle{plain}
\title[Random matrix model for free Meixner laws]
{Random matrix model for free Meixner laws}
\author[R. Lenczewski]{Romuald Lenczewski}
\address{Romuald Lenczewski, \newline
Instytut Matematyki i Informatyki, Politechnika Wroc\l{}awska, \newline
Wybrze\.{z}e Wyspia\'{n}skiego 27, 50-370 Wroc{\l}aw, Poland  \vspace{10pt}}
\email{Romuald.Lenczewski@pwr.wroc.pl}
\begin{abstract}
Applying the concept of matricial freeness which generalizes freeness in free probability
we have recently studied asymptotic joint distributions of symmetric blocks of Gaussian 
random matrices (Gaussian Symmetric Block Ensemble). This approach gives a block refinement of the 
fundamental result of Voiculescu on asymptotic freeness of independent Gaussian random matrices.
In this paper, we show that this framework is natural for constructing a random matrix model for {\it free Meixner laws}. We also demonstrate that 
the ensemble of independent matrices of this type is asymptotically conditionally free with 
respect to the pair of partial traces.
\end{abstract}

\maketitle
\section{Introduction}
It is well-known that free probability is an effective tool in the study of random matrices and their asymptotics.
This approach was originated by Voiculescu in his fundamental paper [17], where he showed that independent Gaussian 
random matrices are asymptotically free (generalized to non-Gaussian entries by Dykema [9]). His result showed that the semicircle law obtained by Wigner [19] as the limit distribution of certain symmetric random matrices can now be viewed as an element of a much more general probability theory involving operator algebras [16].

If the complex-valued Gaussian variables which are entries of the considered random matrices are not 
identically distributed, one has to apply a more general scheme to study their asymptotics.  
One approach is to use operator-valued states and the associated notion of freeness with amalgamation, 
as in the paper of Shlakhtyenko on Gaussian band matrices [15]. This approach was further developed by 
Benayach-Georges [5] who described the asymptotics of blocks of random matrices and introduced
a related additive convolution. Recently, we studied asymptotic joint distributions of symmetric blocks of random 
matrices by means of operatorial methods on Hilbert spaces. For this purpose, we employed a scheme based on arrays of 
scalar-valued states and the associated concept of {\it matricial freeness} introduced in [10].

In particular, we showed in [11,12] that the symmetric blocks of an ensemble of $n\times n$ complex 
Hermitian Gaussian random matrices $Y(u,n)$ with block-identical variances converge in moments under normailzed {\it partial traces} to the 
mixed moments of {\it symmetrized Gaussian operators}, namely
$$
T_{p,q}(u,n)\rightarrow \widehat{\omega}_{p,q}(u)
$$
where $u\in \mathpzc{U}$ and $1\leq p\leq q\leq r$, with $\mathpzc{U}$ being 
an index set enumerating independent matrices. The operators 
$\widehat{\omega}_{p,q}(u)$ are natural symmetrizations of square arrays 
$(\omega_{p,q}(u))$ of {\it matricially free Gaussian operators} 
playing the role of basic Gaussian operators. 
By a partial trace we understand a normalized trace over the subset of basis vectors
related to diagonal blocks. 

In the random matrix context, the corresponding framework is thus a block 
refinement of that used by Voiculescu and is closely related to his idea of decomposition of 
Gaussian random matrices leading to semicircular and circular 
systems [18]. We studied a deformation of this decomposition based on allowing the Gaussian 
variables to have {\it block-identical variances} rather than identical 
and then computing their mixed moments under (normalized) partial traces rather than under the (normalized) complete trace.
We would also like to remark that some results obtained by our methods
can perhaps be suitably reformulated in terms of freeness with amalgamation.

The key parameters of the block refinement are given by $r\times r$ symmetric variance matrices
$V(u)=(v_{p,q}(u))$ associated with symmetric blocks of the matrices $Y(u,n)$, which, in turn are defined
by the partition of the set 
$$
[n]=N_1\cup N_2\cup\ldots \cup N_r
$$ 
into $r$ disjoint intervals (they depend on $n$, but this is supressed in the notation), 
and by the dimension matrix 
$$
D={\rm diag}(d_1, d_2, \ldots , d_r),
$$
whose entries are given by non-negative numbers 
$$
d_j=\lim_{n\rightarrow \infty}\frac{|N_j|}{n}
$$
called {\it asymptotic dimensions}. An important assumption is that 
we allow some of these dimensions to vanish. Note that
in our first paper, where we presented the block model [11],
we assumed that all asymptotic dimensions are positive.

It follows from the asymptotics of random symmetric blocks that the parameters of random matrices 
are encoded in the products of the dimension matrix and the variance matrices, namely
$$
B(u)=DV(u)
$$
and these matrices provide constants associated with blocks of colored non-crossing pair partitions underlying 
the combinatorics of mixed moments of symmetrized Gaussian operators.
Let us add that we take the same dimension matrix for all random matrices.

In comparison with freeness of free probability, matricial freeness gives 
more flexibility in treating such problems of random matrix theory as
\begin{enumerate}
\item
evaluating limit distributions of random matrices,
\item
studying asymptotic properties of random matrix ensembles,
\item
constructing random matrix models for given probability measures,
\end{enumerate}
and in that respect it reminds freeness with amalgamation. Some advantage of our approach 
is that we rely on operators living in Hilbert spaces. This seems quite intutive
especially since computations involve operators which remind free creation and annihilation operators and therefore
their moments can be easily expressed in terms of non-crossing (pair) partitions. A sample of such computations 
is contained in this paper.

In particular, this flexibility allows us to treat sums and products of rectangular random matrices in a unified manner, 
including {\it Wishart matrices} [20] as well as more general products like those leading to {\it free Bessel laws} [4] and 
free products of Marchenko-Pastur [14] distributions with arbitrary shape parameters. In fact, we were able to
compute the moments of the latter in the explicit form (known only in very special cases before) and introduce 
polynomials which can be viewed as {\it multivariate Narayana polynomials} [13].
A number of other new applications to the random matrix theory can be given. In particular, the
matricially free Gaussian operators turned out to be effective in the construction of 
random matrix models for boolean independence, monotone independence and s-freeness [12]. In this paper, we also use 
these operators to construct a simple random matrix model for an important class of probability measures on the real 
line called {\it free Meixner laws} and prove the asymptotic conditional independence of the associated ensembles
of random matrices.

Free Meixner systems of polynomials and the associated family of functionals 
were introduced and studied by Anshelevich [2,3].
Let us remark that free Meixner laws are free analogs 
of {\it classical Meixner laws}. In particular, up to affine transformations, they 
belong to one of the following six classes: free Gaussian (Wigner semicircle), free Poisson (Marchenko-Pastur),
free negative binomial (free Pascal), free Gamma, free binomial and free hyperbolic secant, 
following the terminology of Anshelevich.
Free Meixner laws turn out to display similar properties with respect to free independence as 
do the classical Meixner laws with respect to classical independence 
as Bryc and Bo\.{z}ejko showed in their study of the regression problem [6].

Random matrix models for certain special free Meixner laws are well-known, like the Gaussian Unitary Ensemble 
for the semicircle law, the Wishart Ensemble for the Marchenko-Pastur law or the Jacobi Ensemble for the
free binomial law (see, for instance, [8,17,20]). However, a natural model for the whole class of free Meixner laws has not been 
given in the literature.

The paper is organized as follows. In Section 2, we recall a combinatorial formula for the moments of free Meixner laws.
An operatorial realization of their moments in terms of matricially free Gaussian operators is proved in Section 3. 
A random matrix model for free Meixner laws is
constructed in Section 4. An ensemble of independent random matrices of this type, called the Free Meixner Ensemble, 
is considered in Section 5, where we prove its asymptotic conditional freeness.

\section{Moments of free Meixner laws}
It is well-known that every probability measure on the real line with finite moments of all orders is
characterized by two sequences of Jacobi parameters 
$$
\alpha=(\alpha_{1}, \alpha_{2}, \ldots)\;\;{\rm and}\;\;
\beta =(\beta_{1}, \beta_{2}, \ldots ),
$$
where $\alpha_{n}\in {\mathbb R}$ and
$\beta_{n}\geq 0$ for all $n\in {\mathbb N}\cup\{0\}$, with the condition
that if $\beta_k=0$ for some $k$, then $\beta_m=0$ for all $m>k$.
We will call them {\it Jacobi sequences} and we will 
use the notation $J(\mu)= (\alpha ,\omega)$. The Cauchy transform of $\mu$ can then
be expressed as a continued fraction of the form
$$
G_{\mu}(z)=\cfrac{1}{z-\alpha_{1}-\cfrac{\beta_1}
{z-\alpha_{2}-\cfrac{\beta_{2}}{z-\alpha_{3}-\cfrac{\beta_{3}} {\ldots}}}}
$$
and it is understood that if $\beta_{m}=0$ for some $m$, then the fraction
terminates and, for convenience, we set $\beta_{n}=\alpha_{n}=0$ for all $n>m$.

This continued fraction representation of Cauchy transforms turns out useful 
in our approach. Thus, let us first remark that the family of free Meixner laws 
is the family of probability measures on the real line
associated with the pair of Jacobi sequences of the form 
$$
\alpha=(\alpha_{1}, \alpha_{2}, \alpha_{2}, \ldots )\;\;{\rm and}\;\;
\beta=(\beta_{1}, \beta_{2}, \beta_{2}, \ldots ),
$$
i.e. they are constant starting from the second level of the corresponding continued fractions.
If a free Meixner law corresponds to the pair of Jacobi sequences of the above form, we will say that
it corresponds to $(\alpha_1,\alpha_2, \beta_1, \beta_2)$.
In particular, if $\alpha_{1}=0$ and $\beta_{1}=1$, we obtain the {\it standard} free Meixner laws 
with mean zero and variance one. In that case, the absolutely continuous part of the associated measure $\mu$
takes the form
$$
d\mu(x)=\frac{\sqrt{4\beta_2-(x-\alpha_2)^{2}}}{2\pi(\beta_2-1)x^2+\alpha_2x+1}
$$
on $[\alpha_2-2\sqrt{\beta_2},\alpha_2+2\sqrt{\beta_2}]$, the measure can also have one or two atoms.

There is a useful combinatorial formula which expresses moments of probability measures on the real line 
in terms of non-crossing partitions consisting of 1-blocks (singletons) and 2-blocks (pairs).
Namely, let $\mathcal{NC}^{1,2}_m$ be the set of non-crossing partitions of the set 
$[m]=:\{1,2, \ldots , m\}$ consisting of singletons and pairs, namely
$$
\pi=\{\pi_1, \pi_2, \ldots , \pi_k\}\in \mathcal{NC}_m
$$
where each $\pi_j$ contains one or two elements, respectively, 
and it is not possible to have two different 2-blocks $\pi_i=\{p,q\}$ and $\pi_j=\{r,s\}$, for which $p<r<q<s$. 

In any non-crossing partition $\pi\in \mathcal{NC}_{m}$, if we put all numbers from the set $[m]$ in order 
and draw lines connecting all numbers which belong to the same block, the lines corresponding to different blocks cannot intersect each other. 
Further, its block $\pi_i$ is {\it outer} with respect to the block $\pi_j$ if there exist
$r,s\in \pi_i$ such that for each $p\in \pi_j$ it holds that $r<p<s$. If $\pi$ consists of singletons and pairs, it is clear that
any outer block must be a pair. We say that the block $\pi_i$ of $\pi\in \mathcal{NC}$ has {\it depth} 
$\mathpzc{d}(\pi_i)=\mathpzc{d}(i)$ if it has $\mathpzc{d}(i)-1$ outer blocks. 
Thus, blocks which do not have outer blocks are assumed to have depth one. Note that if a block $\pi_i$ has 
at least one outer block, we can choose among them the one which lies immediately above $\pi_i$ and we will call
it its {\it nearest outer block}.

If $\mu$ is a probability measure on the real line with all moments finite and the pair of Jacobi sequences
$J(\mu)=(\alpha, \beta)$, its $n$-th moment is given by the combinatorial formula
$$
M_{n}(\mu)=\sum_{\pi\in \mathcal{NC}^{1,2}_n}\;\;\prod_{i: |\pi_{i}|=1}\alpha_{\mathpzc{d}(i)}\prod_{j:|\pi_{j}|=2}\beta_{\mathpzc{d}(j)},
$$
i.e. each block of depth $\mathpzc{d}$ of every $\pi\in \mathcal{NC}^{1,2}_m$ contributes $\alpha_{\mathpzc{d}}$ or $\beta_{\mathpzc{d}}$
if it is a singleton or a pair, respectively. This formula was first discovered by Cabanal-Duvillard and Ionesco
for symmetric measures [7]. In that case, the first Jacobi sequence $\alpha$ vanishes and 
only pair partitions appear in the formula. The general version is due to Accardi and Bo\.{z}ejko [1].

\section{Operatorial realization}
We will use matricially free Gaussian operators living in the {\it matricially free Fock space of tracial type} 
introduced in [12] to find a realization of moments of free Meixner laws. 
This Fock space is a generalization of the matricially free Fock space $\mathcal{M}$ introduced in [11]. 

For the purposes of this article, it suffices to consider the special case when 
$$
{\mathcal M}={\mathcal M}_{1}\oplus {\mathcal M}_{2},
$$
where both ${\mathcal M}_{1}$ and ${\mathcal M}_{2}$ are Hilbert space direct sums
\begin{eqnarray*}
{\mathcal M}_{1}&=&{\mathbb C}\Omega_1\oplus 
\bigoplus_{k=0}^{\infty}({\mathcal H}_{2}^{\otimes k}\otimes {\mathcal H}_{1}),\\
{\mathcal M}_{2}&=&{\mathbb C}\Omega_2\oplus 
\bigoplus_{k=1}^{\infty}{\mathcal H}_{2}^{\otimes k},\\
\end{eqnarray*}
where $\Omega_1, \Omega_2$ are unit vectors, 
${\mathcal H}_{j}={\mathbb C}e_{j}$ for $j\in \{1,2\}$, where
$e_1,e_2$ are unit vectors, and ${\mathcal H}^{\otimes 0}\otimes {\mathcal H}_{1}={\mathcal H}_{1}$.
The space ${\mathcal M}$ is endowed with the canonical inner product. 

Using the canonical basis of this Fock space,
$$
{\mathpzc B}= \{\Omega_1, \Omega_2, 
e_{2}^{\otimes k}\otimes e_{1}, \;e_{2}^{\otimes l}:\, k\in {\mathbb N}\cup\{0\}, l\in {\mathbb N}\},
$$
we define creation operators $\wp_{1},\wp_{2}\in B({\mathcal M})$ as follows. Let $(\beta_1, \beta_2)$ 
be a pair of nonnegative numbers. We set
$$
\wp_{1}\Omega_1=\sqrt{\beta_{1}}\,e_{1},
$$
and we assume that $\wp_{1}$ sends the remaining basis vectors to zero. 
In turn, $\wp_2$ sends $\Omega_1$ to zero and otherwise,
\begin{eqnarray*}
\wp_{2}\Omega_2&=&\sqrt{\beta_{2}}\,e_{2}\\
\wp_{2}(e_{2}^{\otimes k})&=&\sqrt{\beta_{2}}\,e_{2}^{\otimes (k+1)}\\
\wp_{2}(e_{2}^{\otimes l}\otimes e_{1})&=&\sqrt{\beta_{2}}\,(e_{2}^{\otimes (l+1)}\otimes e_{1})
\end{eqnarray*}
for any $k\in {\mathbb N}, l\in {\mathbb N}\cup \{0\}$.
By $\wp_{1}^{*}$ and $\wp_{2}^{*}$ we denote 
the adjoints of $\wp_1$ and $\wp_2$, respectively, and sums of the form
$$
\omega_{1}=\wp_{1}+\wp_{1}^{*}\;\;{\rm and}\;\;\omega_{2}=\wp_{2}+\wp_{2}^{*}
$$
will be the corresponding Gaussian operators. Note that ${\mathcal M}_{j}$ is invariant 
with respect to $\wp_{i}, \wp_{i}^{*}, \omega_{i}$ for any $i,j\in \{1,2\}$.

In particular, if we set $\beta_1=\beta_2=1$, then the restrictions
$$
(\wp_1+\wp_2)|{\mathcal M}_{1}\;\;{\rm and}\;\;\wp_2|{\mathcal M}_{2}
$$ 
can be identified with the standard free creation operators living in ${\mathcal M}_{1}$ and ${\mathcal M}_{2}$, respectively, 
and both spaces are isomorphic to the free Fock space over the one-dimensional Hilbert space.

\begin{Remark}
{\rm 
Our Fock space ${\mathcal M}$ is a special case of the 
{\it matricially free Fock space of tracial type} associated with an array 
$({\mathcal H}_{p,q})$ of Hilbert spaces, by which we understand the Hilbert space direct sum
$$
{\mathcal M}= \bigoplus_{q=1}^{r} {\mathcal M}_{q},
$$
where each summand is of the form
$$
{\mathcal M}_{q}={\mathbb C}\Omega_{q}\oplus \bigoplus_{m=1}^{\infty}
\bigoplus_
{(p_1,p_2)\neq (p_2,p_3)\neq \ldots \neq (p_m,q)}
{\mathcal F}_{p_1,p_2}^{\,0}\otimes {\mathcal F}_{p_2,p_3}^{\,0}
\otimes \ldots \otimes 
{\mathcal F}_{p_m,q}^{\,0}
$$
with tensor products built from free and boolean Fock spaces 
$$
{\mathcal F}_{p,q}^{0}=\left\{\begin{array}{ll}
\bigoplus_{k=1}^{\infty}{\mathcal H}_{q,q}^{\otimes k} & {\rm if}\;p=q\\
{\mathcal H}_{p,q} & {\rm if}\;p\neq q 
\end{array}
\right.,
$$
with vacuum spaces subtracted.
In this paper, we suppose the array $({\mathcal H}_{p,q})$ consists of only two one-dimensional Hilbert spaces
${\mathcal H}_{2,1}={\mathcal H}_1$ and ${\mathcal H}_{2,2}={\mathcal H}_{2}$.
Clearly, an assymmetry in $({\mathcal H}_{p,q})$ leads to an assymetry in the 
definitions of ${\mathcal M}_{1}$ and ${\mathcal M}_{2}$.}
\end{Remark}

\begin{Remark}
{\rm We can identify the creation operators $\wp_1, \wp_2$ with 
the {\it matricially free creation operators}
$$
\wp_{1}=\wp_{2,1}\;\;{\rm and}\;\;\wp_{2}=\wp_{2,2},
$$
where we use the matricial two-index notation of [10,11]. This notation is often helpful 
(and will be used when we refer to the results of these papers)
since the second index shows onto which basis vectors the operators act non-trivially (it must match the first index of
the basis vector). Therefore,
$\wp_{p,q}$ acts non-trivially only onto $\Omega_q$ and tensor products which begin with $e_{q,r}$
for any $r$. Thus, for instance,
$$
\wp_{2,1}\Omega_2=0,\; \wp_{2,1}e_{2,1}=0, \;\wp_{2,1}(e_{2,2}\otimes e_{2,1})=0, \;\wp_{2,2}\Omega_1=0, 
$$
which stands behind the definition of $\wp_1, \wp_2$ (the fact that $\wp_{1,2}$ and $\wp_{1,1}$ are not 
used makes the one-index notation feasible). We also have $\wp_{1}^{*}=\wp_{2,1}^{*}$, $\wp_{2}^{*}=\wp_{2,2}^{*}$, 
with the corresponding scalars $\beta_1=b_{2,1}$ and $\beta_2=b_{2,2}$. 
In turn, $\omega_1$ and $\omega_2$ can be identified with the corresponding 
{\it matricially free Gaussian operators} $\omega_{2,1}$ and $\omega_{2,2}$, respectively.
Details on the arrays of such operators can be found in [10,11].}
\end{Remark}

Using these operators, we can define operators in $B({\mathcal M}_1)$ whose distributions
in the state $\Psi_1$ defined by the vector $\Omega_1$ are free Meixner laws.
For that purpose, the subspace ${\mathcal M}_{2}$ is not needed yet.

\begin{Theorem}
If $\mu$ is the free Meixner law corresponding to $(\alpha_1, \alpha_2,\beta_1,\beta_2)$, where $\beta_1\neq 0$ and 
$\beta_2 \neq 0$, then its $m$-th moment is given by
$$
M_{m}(\mu)=\Psi_1((\omega+\gamma)^{m}),
$$
where 
$$
\omega=\omega_{1}+\omega_{2}
$$
and 
$$
\gamma=(\alpha_{2}-\alpha_1)(\beta_{1}^{-1}\wp_{1}\wp_{1}^{*}+
\beta_{2}^{-1}\wp_{2}\wp_{2}^{*})+\alpha_{1},
$$
and $\Psi_1$ is the state defined by the vector $\Omega_1$.
\end{Theorem}
{\it Proof.}
Let us first analyze the moments of $\omega$ since these were studied in [12]
in the general case of matricially free Gaussian operators.
The operator $\omega$ can be identified with 
$$
\omega=\omega_{2,1}+\omega_{2,2}
$$ 
by Remark 3.2. Of course, if we set $\omega_{1,2}=\omega_{1,1}=0$, we can use the
combinatorial formula for the moments of the {\it Gaussian pseudomatrix},
$$
\omega=\sum_{1\leq p,q\leq 2}\omega_{p,q}
$$
associated with a $2\times 2$ array $(\omega_{p,q})$, in which we express these moments in terms of 
colored non-crossing pair partitions [11, Lemma 4.1].

By a colored non-crossing pair partition we shall understand a pair $(\pi, f)$, where $\pi=\{\pi_1,\pi_2, \ldots , \pi_s\}$ 
is a non-crossing pair partition and $f$ is a function on the set of its blocks with values in the set $[r]$. If we
draw an additional 2-block which is outer with  respect to all blocks of $\pi$, called the {\it imaginary block}, 
and we color it by $q$, we obtain the set of colored non-crossing pair partitions 
$\mathcal{NC}_{m,q}^{\,2}[r]$ colored by $[r]$ under condition that the imaginary block is 
colored by $q$. Then, we have 
$$
\Psi_{q}(\omega^m)=\sum_{(\pi,f)\in \mathcal{NC}_{m,q}^{\,2}[r]}b_{q}(\pi,f),
$$
where the summation is over the empty set if $m$ is odd, and
$$
b_q(\pi,f)=b_{q}(\pi_1,f)b_{q}(\pi_2,f)\ldots b_{q}(\pi_s,f)
$$
if $m=2s$, where $\pi=\{\pi_1, \pi_2,\ldots, \pi_s\}$ and 
$$
b_{q}(\pi_k,f)=b_{i,j}
$$ 
whenever block $\pi_k$ is colored by $i$ and its nearest outer block is colored by $j$.
In this formulation, we set $b_{1,1}=b_{1,2}=0$ since there is no $\omega_{1,1}$ or $\omega_{1,2}$, 
but formally it holds for all colorings.

If we set $q=1$, which refers to our theorem, the imaginary block gets 
colored by $1$. Moreover, since $b_{1,2}=b_{1,1}=0$, the non-vanishing contribution to
$\Psi_1(\omega^{m})$ comes only from those colored partitions $(\pi,f)\in \mathcal{NC}_{m,q}^{\,2}[r]$ 
in which each block of $\pi$ is colored by $2$. In fact, if some block $\pi_k$ 
was colored by $1$ and its nearest outer block (including the imaginary block) was colored by $1$ or $2$, then 
the corresponding $b(\pi_k,f)$ would have to be $b_{1,1}$ or $b_{1,2}$, but these vanish. This means that
to each block of depth one we assign the number $b_{2,1}=\beta_{1}$ since the imaginary block is colored by $1$
and it is its nearest outer block, whereas to each block of depth greater than one 
we assign the number $b_{2,2}=\beta_2$ since each block of $\pi$ is colored by $2$.
Namely
$$
b_{1}(\pi_k,f)=\left\{\begin{array}{ll}
\beta_{1} & {\rm if}\; \mathpzc{d}(\pi_k)=1\\
\beta_{2} & {\rm if}\; \mathpzc{d}(\pi_k)>1
\end{array}
\right..
$$
This gives
$$
\Psi_1((\omega_{2,2}+\omega_{2,1})^{m})=
\sum_{\pi \in \mathcal{NC}_{m}^{2}}b_{1}(\pi,f)
$$
since in this case the set $\mathcal{NC}_{m,1}^{2}[2]$ of colored non-crossing pair
partitions of $[m]$ with the imaginary block colored by $1$ reduces to 
the set in which all blocks colored by $2$, which is in bijection with 
$\mathcal{NC}_{m}^{2}$. Switching back to the notations of this paper, we thus have
$$
\Psi_1(\omega^{m})=
\sum_{\pi \in \mathcal{NC}_{m}^{2}}\beta_{1}^{|B_{1}(\pi)|}\beta_{2}^{|B_{2}(\pi)|},
$$
where $B_{1}(\pi)$ and $B_{2}(\pi)$ are the sets of 2-blocks of $\pi$ of depth $1$ nad of depth greater than $1$, respectively.

In fact, the above formula for the moments of $\omega$ can be proved directly without invoking the general statement of [11, Lemma 4.1]. It suffices to observe that $(\pi,f)$ is uniquely determined by the sequence $\epsilon=(\epsilon_1,\epsilon_2, \ldots , \epsilon_m)$
which appears in nonvanishing mixed moments of creation and annihilation operators of type
$$
\Psi_1(\wp_{q_1}^{\epsilon_1}\wp_{q_2}^{\epsilon_2}\ldots \wp_{q_m}^{\epsilon_m}),
$$
where $\epsilon_k\in \{1,*\}$ since the choice of $\epsilon$ 
uniquely determines the tuple $(q_1, q_2, \ldots , q_m)$ due to 
the 0-1 action of $\wp_{1}$ and $\wp_{2}$ and their adjoints. 
Namely, only $\wp_1$ acts non-trivially onto $\Omega_1$, giving $e_{1}$, 
which corresponds to the right leg of each block of depth $1$ 
(its adjoint corresponds to its left leg since it sends $e_{1}$ into $\Omega_1$). 
In turn, $\wp_{2}$ acts non-trivially onto each basis element of $\mathpzc{B}$ except $\Omega$ 
and thus it corresponds to the right leg of each block of depth greater than $1$ 
(its adjoint corresponds to its left leg). Therefore, each block of $\pi$ of depth $1$
is associated with the pair $(\wp_{1}^{*},\wp_{1})$ producing $\beta_1$, whereas the remaining blocks
are associated with the pair $(\wp_{2}^{*}, \wp_{2})$ producing $\beta_2$.

It remains to check what happens when we replace $\omega$
by $\omega+\gamma$. Observe that 
$$
\wp_{1}\wp_{1}^{*}=\beta_{1}P_{1}\;\;{\rm and}\;\;\wp_{2}\wp_{2}^{*}=\beta_{2}P_{2},
$$
where $P_1$ is the canonical projection onto ${\mathcal H}_{1}$ and $P_2$ is the 
canonical projection onto the subspace
$$
{\mathcal F}_{2}=\bigoplus_{k=1}^{\infty}({\mathcal H}_{2}^{\otimes k}\otimes {\mathcal H}_{1}),
$$
respectively. Therefore,
$$
\gamma=\alpha_1P+\alpha_2(P_1+P_2),
$$
where $P$ is the canonical projection onto ${\mathbb C}\Omega_1$,
which means that $\gamma$ is diagonal in the basis $\mathpzc{B}$, 
namely it multiplies $\Omega_1$ and all vectors from $\mathpzc{B}\setminus {\Omega}_1$ by $\alpha_1$ and $\alpha_2$, respectively.
Therefore, if we are given a mixed moment
$$
\Psi_1(\wp_{q_1}^{\epsilon_1}\wp_{q_2}^{\epsilon_2}\ldots \wp_{q_k}^{\epsilon_k}),
$$
associated with a non-crossing pair partition $\pi$ of the set $[k]$, each mixed moment of the form
$$
\Psi_1(\gamma^{n_0}\wp_{q_1}^{\epsilon_1}\gamma^{n_1}\wp_{q_2}^{\epsilon_2}\gamma^{n_2}\ldots \wp_{q_k}^{\epsilon_k}\gamma^{n_k}),
$$
where $n_0,n_1, \ldots , n_k$ are non-negative integers such that 
$$
k+n_0+n_1+\ldots +n_k=m,
$$
which appears when we compute the $m$-th moment of $\omega+\gamma$, 
is naturally associated with a non-crossing partition $\widetilde{\pi}$ 
of the set $[m]$ obtained from $\pi$ by adding $m-k$ singletons in such a way that
$n_j$ singletons are placed right after the number $j$, with 
$n_0$ singletons placed before the number 1 belonging to the first pair.
In this fashion we obtain all non-crossing partitions of $[m]$ which have 
$m-k$ singletons and $k$ pairs. Further, each $\widetilde{\pi}\in\mathcal{NC}_{m}^{1,2}$ 
is obtained exactly once in this fashion from some $\pi\in \mathcal{NC}_{k}^{2}$.

Moreover, to each singleton of depth $1$ 
we assign $\alpha_1$ and to each singleton of depth greater than $1$ we assign $\alpha_2$ in view of the 
diagonal form of $\gamma$ in the basis $\mathpzc{B}$. 
Therefore, we obtain
$$
\Psi_1((\omega+\gamma)^{m})=\sum_{\pi\in \mathcal{NC}_{m}^{1,2}}\;\;
\alpha_{1}^{|S_{1}(\pi)|}\alpha_{2}^{|S_{2}(\pi)|}\beta_{1}^{|B_{1}(\pi)|}\beta_{2}^{|B_{2}(\pi)|},
$$
where $S_1(\pi)$ and $S_{2}(\pi)$ are the sets of singletons of depth $1$ and of depth 
greater than $1$ in $\pi$, respectively.
As we know from the combinatorial formula for the moments given in the Introduction, 
this is the $m$-th moment of the free Meixner law. This completes the proof.\hfill $\blacksquare$\\

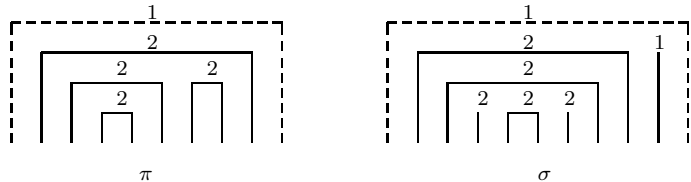
\begin{figure}
\unitlength=1mm
\special{em:linewidth 0.4pt}
\linethickness{0.4pt}
\begin{picture}(120.00,30.00)(5.00,5.00)
\put(22.00,10.00){\line(0,1){1.25}}
\put(58.00,10.00){\line(0,1){1.25}}
\put(22.00,12.00){\line(0,1){1.25}}
\put(58.00,12.00){\line(0,1){1.25}}
\put(22.00,14.00){\line(0,1){1.25}}
\put(58.00,14.00){\line(0,1){1.25}}
\put(22.00,16.00){\line(0,1){1.25}}
\put(58.00,16.00){\line(0,1){1.25}}
\put(22.00,18.00){\line(0,1){1.25}}
\put(58.00,18.00){\line(0,1){1.25}}
\put(22.00,20.00){\line(0,1){1.25}}
\put(58.00,20.00){\line(0,1){1.25}}
\put(22.00,22.25){\line(0,1){1.25}}
\put(58.00,22.25){\line(0,1){1.25}}
\put(22.00,24.50){\line(0,1){1.50}}
\put(58.00,24.50){\line(0,1){1.50}}

\put(26.00,10.00){\line(0,1){12.00}}
\put(30.00,10.00){\line(0,1){8.00}}
\put(34.00,10.00){\line(0,1){4.00}}
\put(38.00,10.00){\line(0,1){4.00}}
\put(42.00,10.00){\line(0,1){8.00}}
\put(46.00,10.00){\line(0,1){8.00}}
\put(50.00,10.00){\line(0,1){8.00}}
\put(54.00,10.00){\line(0,1){12.00}}

\put(40.00,22.50){$\scriptstyle{2}$}
\put(36.00,19.00){$\scriptstyle{2}$}
\put(36.00,15.00){$\scriptstyle{2}$}
\put(48.00,19.00){$\scriptstyle{2}$}
\put(40.00,26.50){$\scriptstyle{1}$}

\put(22.00,26.00){\line(1,0){1.25}}
\put(24.00,26.00){\line(1,0){1.25}}
\put(26.00,26.00){\line(1,0){1.25}}
\put(28.06,26.00){\line(1,0){1.25}}
\put(30.12,26.00){\line(1,0){1.25}}
\put(32.18,26.00){\line(1,0){1.25}}
\put(34.26,26.00){\line(1,0){1.25}}
\put(36.32,26.00){\line(1,0){1.25}}
\put(38.40,26.00){\line(1,0){1.25}}
\put(40.48,26.00){\line(1,0){1.25}}
\put(42.56,26.00){\line(1,0){1.25}}
\put(44.62,26.00){\line(1,0){1.25}}
\put(46.70,26.00){\line(1,0){1.25}}
\put(48.70,26.00){\line(1,0){1.25}}
\put(50.72,26.00){\line(1,0){1.25}}
\put(52.75,26.00){\line(1,0){1.25}}
\put(54.75,26.00){\line(1,0){1.25}}
\put(56.75,26.00){\line(1,0){1.25}}

\put(26.00,22.00){\line(1,0){28.00}}
\put(30.00,18.00){\line(1,0){12.00}}
\put(34.00,14.00){\line(1,0){4.00}}
\put(46.00,18.00){\line(1,0){4.00}}
\put(39.00,5.00){$\scriptstyle{\pi}$}
\put(72.00,10.00){\line(0,1){1.25}}
\put(112.00,10.00){\line(0,1){1.25}}
\put(72.00,12.00){\line(0,1){1.25}}
\put(112.00,12.00){\line(0,1){1.25}}
\put(72.00,14.00){\line(0,1){1.25}}
\put(112.00,14.00){\line(0,1){1.25}}
\put(72.00,16.00){\line(0,1){1.25}}
\put(112.00,16.00){\line(0,1){1.25}}
\put(72.00,18.00){\line(0,1){1.25}}
\put(112.00,18.00){\line(0,1){1.25}}
\put(72.00,20.00){\line(0,1){1.25}}
\put(112.00,20.00){\line(0,1){1.25}}
\put(72.00,22.25){\line(0,1){1.25}}
\put(112.00,22.25){\line(0,1){1.25}}
\put(72.00,24.50){\line(0,1){1.50}}
\put(112.00,24.50){\line(0,1){1.50}}

\put(76.00,10.00){\line(0,1){12.00}}
\put(80.00,10.00){\line(0,1){8.00}}
\put(84.00,10.00){\line(0,1){4.00}}
\put(88.00,10.00){\line(0,1){4.00}}
\put(92.00,10.00){\line(0,1){4.00}}
\put(96.00,10.00){\line(0,1){4.00}}
\put(100.00,10.00){\line(0,1){8.00}}
\put(104.00,10.00){\line(0,1){12.00}}
\put(108.00,10.00){\line(0,1){12.00}}

\put(90.00,22.50){$\scriptstyle{2}$}
\put(90.00,19.00){$\scriptstyle{2}$}
\put(90.00,15.00){$\scriptstyle{2}$}
\put(90.00,26.50){$\scriptstyle{1}$}
\put(84.00,15.00){$\scriptstyle{2}$}
\put(95.50,15.00){$\scriptstyle{2}$}
\put(107.50,22.50){$\scriptstyle{1}$}

\put(72.00,26.00){\line(1,0){1.25}}
\put(74.00,26.00){\line(1,0){1.25}}
\put(76.00,26.00){\line(1,0){1.25}}
\put(78.06,26.00){\line(1,0){1.25}}
\put(80.12,26.00){\line(1,0){1.25}}
\put(82.18,26.00){\line(1,0){1.25}}
\put(84.26,26.00){\line(1,0){1.25}}
\put(86.32,26.00){\line(1,0){1.25}}
\put(88.40,26.00){\line(1,0){1.25}}
\put(90.48,26.00){\line(1,0){1.25}}
\put(92.56,26.00){\line(1,0){1.25}}
\put(94.62,26.00){\line(1,0){1.25}}
\put(96.70,26.00){\line(1,0){1.25}}
\put(98.70,26.00){\line(1,0){1.25}}
\put(100.72,26.00){\line(1,0){1.25}}
\put(102.75,26.00){\line(1,0){1.25}}
\put(104.75,26.00){\line(1,0){1.25}}
\put(106.75,26.00){\line(1,0){1.25}}
\put(108.75,26.00){\line(1,0){1.25}}
\put(110.75,26.00){\line(1,0){1.25}}

\put(76.00,22.00){\line(1,0){28.00}}
\put(80.00,18.00){\line(1,0){20.00}}
\put(88.00,14.00){\line(1,0){4.00}}
\put(96.00,18.00){\line(1,0){4.00}}

\put(92.00,5.00){$\scriptstyle{\sigma}$}

\end{picture}
\caption{Examples of colored non-crossing partitions.}
\end{figure}

\begin{Example}
{\rm Let us give some examples of non-crossing partitions and the associated mixed moments. 
The diagrams are given in Figure 1. The partition $\pi$ consists of 4 pairs, namely
$\pi_1=\{1,8\}$, $\pi_2=\{2,5\}$, $\pi_3=\{3,4\}$, $\pi_4=\{6,7\}$, with the imaginary block 
marked with a dotted line. There exists exactly one mixed moment of creation and annihilation operators that corresponds to this
partition, namely we must have $\epsilon=(*,*,*,1,1,*,1,1)$ and the corresponding moment (the only non-trivial one which corresponds to
this $\epsilon$) is 
$$
\Psi_1(\wp_{1}^{*}\wp_{2}^{*}\wp_{2}^{*}\wp_{2}\wp_{2}\wp_{2}^{*}\wp_{2}\wp_{1})=\beta_1 \beta_2^{3}
$$ 
since $\wp_{1}$ is the only creation operator which acts non-trivially onto $\Omega_1$, giving $e_{1}$,
and $\wp_{2}$ is the only creation operator which acts non-trivially onto $e_{1}$ and $e_2\otimes e_1$, 
giving $e_{2}\otimes e_{1}$ and $e_2^{\otimes 2}\otimes e_1$, respectively. 
Next, $\wp_{1}^{*}$ is the only annihilation operator which acts non-trivially onto $e_1$, 
whereas $\wp_2^{*}$ is the only annihilation operator which acts non-trivially 
onto $e_2\otimes e_1$ and $e_2^{\otimes 2}\otimes e_1$.

The partition $\sigma$ contains 3 pairs and 3 singletons, namely 
$\sigma_1=\{1,8\}$, $\sigma_2=\{2,7\}$, 
$\sigma_3=\{3\}$, $\sigma_4=\{4,5\}$, $\sigma_5=\{6\}$, $\sigma_{6}=\{9\}$. 
We assign the color $1$ to all singletons of depth one and the color $2$ to all 
remaining singletons. The colors assigned to singletons are 
to some extent arbitrary (they did not appear in [11,12], where we considered pair partitions only), 
but it is convenient to color all singletons of depth $1$ by $1$ and the remaining ones by $2$
since this corresponds to the right Jacobi coefficients. The associated mixed moment is 
$$
\Psi_1(\wp_{1}^{*}\wp_{2}^{*}\gamma\wp_{2}^{*}\wp_{2}\gamma\wp_2\wp_{1}\gamma)=\alpha_1\alpha_2^2\beta_1\beta_2^2,
$$
where the 2-blocks are associated with the pairs $(\wp_{1}^{*}, \wp_1)$ and $(\wp_2^{*},\wp_2)$, which 
produce $\beta_1$ and $\beta_2$, respectively (like in the case of $\pi$), 
whereas the singletons are associated with $\gamma$, which produces $\alpha_1$ in the case of 
$\{9\}$ (since in this case $\gamma$ acts onto $\Omega_1$), and
$\alpha_2$ in the case of $\{3\}$ and $\{9\}$ (since in this case $\gamma$ acts onto $e_2\otimes e_1$).
}
\end{Example}

If $\beta_1=\beta_2=0$, we set $\omega_1=0$ and $\gamma_1=\alpha$ which leads to
the Dirac measure at $\alpha_1$. In turn, the case $\beta_2=0$ is treated below. 
\begin{Corollary}
If $\mu$ is the free Meixner law corresponding to $(\alpha_1, \alpha_2,\beta_1,0)$,
then its $m$-th moment is given by
$$
M_{m}(\mu)=\Psi_1((\omega_1+\gamma_1)^{m}),
$$
where 
$$
\gamma_1=(\alpha_{2}-\alpha_1)\beta_{1}^{-1}\wp_{1}\wp_{1}^{*}+\alpha_{1}
$$
and $\Psi_1$ is the state defined by the vector $\Omega_1$. 
\end{Corollary}
{\it Proof.}
It suffices to observe that if we disregard $\wp_2$ and $\wp_2^{*}$ in all 
computations in the proof of Theorem 3.1, then $\beta_2$ disappears from the 
formula for the moments of $\omega+\gamma$ under $\Psi_1$.  \hfill $\blacksquare$\\

Finally, we would like to compute the moments of $\omega+\gamma$ in the state $\Psi_2$.
Observe that $\wp_{2,1}$ vanishes on ${\mathcal M}_{2}$ and therefore this reduces 
to the computation of moments of a slightly simpler operator.
\begin{Corollary}
If $\mu$ is the free Meixner law corresponding to $(\alpha_1, \alpha_2,\beta_2,\beta_2)$,
where $\beta_2>0$, then its $m$-th moment is given by
$$
M_{m}(\mu)=\Psi_2((\omega_2+\gamma_2)^{m}),
$$
where 
$$
\gamma_2=(\alpha_{2}-\alpha_1)\beta_{2}^{-1}\wp_{2}\wp_{2}^{*}+\alpha_{1}
$$
and $\Psi_2$ is the state defined by the vector $\Omega_2$. 
\end{Corollary}
{\it Proof.}
Observe that the action of $\wp_2, \wp_{2}^{*}$ on ${\mathcal M}_{2}$ is exactly the 
same as that of the free creation and annihilation operators, respectively, 
on the free Fock space. This means that the moments of $\omega_2$ under $\Psi_2$ agree with 
the moments of the (centered) semicircle law with variance $\beta$, i.e. each moment 
of even order $m=2s$ is equal to $\beta^s$ times the Catalan number $C_s$. Represent $C_s$ as
the sum over $\mathcal{NC}_{m}^{2}$ and observe that
if we replace $\omega_2$ by $\omega_2+\gamma_2$, the effect is that $\mathcal{NC}_{m}^{2}$ 
gets replaced by $\mathcal{NC}_{m}^{1,2}$ as in the proof of Theorem 3.1, with
singletons of depth $1$ and $2$ contributing $\alpha_1$ and $\alpha_2$, respectively.
This gives the combinatorial formula for the $m$-th moment of the free Meixner law corresponding to 
$(\alpha_1, \alpha_2, \beta_2,\beta_2)$.
\hfill $\blacksquare$

\section{Random matrix model}

Using our results on asymptotic distributions of random symmetric blocks and 
Theorem 3.1, we can now construct a random matrix model for free Meixner laws.

Consider the sequence of Gaussian Hermitian random matrices $Y(n)$, where $n\in {\mathbb N}$, 
under the assumptions of [11, Theorem 5.1]. Namely, we assume that $Y(n)$ is a complex
Gaussian $n\times n$ random matrix of the block form
$$
Y(n)=\left(
\begin{array}{rr}
A(n)& B(n)\\
C(n)& D(n)
\end{array}
\right)
$$
where the off-diagonal blocks are adjoints of each other, whereas the diagonal blocks are Hermitian
and the sizes of blocks are defined by the partition of the set $[n]=\{1,2,\ldots, n\}$,
$$
[n]=N_1\cup N_2, \;\; {\rm where}\;\;N_1\cap N_2=\emptyset
$$
and 
$$
d_{1}=\lim_{n\rightarrow \infty}\frac{N_1}{n}=0\;\;{\rm and}\;\;
d_{2}=\lim_{n\rightarrow \infty}\frac{N_2}{n}=1,
$$
which corresponds to the situation in which 
\begin{enumerate}
\item
the sequence $(D(n))$ is {\it balanced},
\item
the sequence of symmetric blocks built from 
$(B(n))$ and $(C(n))$ is {\it unbalanced},
\item
the sequence $(A(n))$ is {\it evanescent},
\end{enumerate}
according to the natural terminology introduced in [12]. Since $(A(n))$ is evanescent, we can equivalently 
assume that each block of this sequence vanishes.

Using the notation of [12], where blocks are equipped with indices, 
we have 
$$
A(n)=S_{1,1}(n), B(n)=S_{1,2}(n), C(n)=S_{2,1}(n), D(n)=S_{2,2}(n).
$$
It is convenient to identify all blocks $S_{p,q}(n)$ as well as the symmetric blocks
$$
T_{p,q}(n)=\left\{
\begin{array}{ll}
S_{q,q}(n) &{\rm if}\,p=q\\
S_{p,q}(n)+S_{q,p}(n) & {\rm if}\,p\neq q
\end{array}
\right.
$$
with their embeddings in the algebra of $n\times n$ matrices, so that we can decompose
matrices in terms of their blocks, namely
$$
Y(n)=\sum_{p,q}S_{p,q}(n)=\sum_{p\leq q}T_{p,q}(n),
$$
which allows us to write the mixed moments of blocks under any partial trace
$\tau_j(n)$ over basis vectors of ${\mathbb C}^{n}$ indexed by the set $N_j$.

Shortly speaking, we shall assume that the matrices $Y(n)$ are Gaussian Hermitian random 
matrices with block-identically distributed entries. More explicitly, we assume that 
\begin{enumerate}
\item
each entry  $Y_{i,j}(n)$ of $Y_(n)$ is a complex Gaussian random variable
of the form 
$$
Y_{i,j}(n)={\rm Re}Y_{i,j}(n)+i{\rm Im}Y_{i,j}(n),
$$
\item
the family 
$$
\{{\rm Re}Y_{i,j}(n), {\rm Im}Y_{i,j}(n): 1\leq i \leq j\leq n\}
$$
is independent for any $n$,
\item
the real-valued Gaussian variables have mean zero and 
$$
{\mathbb E}(\overline{Y_{i,j}(n)}Y_{i,j}(n))=\frac{v_{p,q}}{n}
$$
whenever $(i,j)\in N_p\times N_q$ for $p,q\in \{1,2\}$, where the variance matrix $V=(v_{p,q})$
is symmetric.
\end{enumerate}

\begin{Theorem}
Under the above assumptions, let $\tau_1(n)$ be the partial normalized trace over the set of 
first $N_1$ basis vectors and let $\beta_1=v_{2,1}>0$ and $\beta_2=v_{2,2}>0$. Then
$$
\lim_{n\rightarrow \infty}\tau_1(n)\left((M(n))^{m}\right)=\Psi_1((\omega+\gamma)^{m})
$$
where 
$$
M(n)=Y(n)+\alpha_1I_1(n)+\alpha_2I_2(n)
$$
for any $n\in {\mathbb N}$, where $I(n)=I_1(n)+I_2(n)$ is the decomposition of the 
$n\times n$ unit matrix induced by the partition $[n]=N_1\cup N_2$ 
and $\omega, \gamma$ are given by Theorem 3.1.
\end{Theorem}
{\it Proof.}
We decompose $Y(n)$ in terms of symmetric random blocks as
$$
Y(n)=T_{1,2}(n)+T_{1,1}(n)+T_{2,2}(n)
$$
and therefore, by [11, Theorem 5.1], the moments of $Y(n)$ under any partial trace, including $\tau_1(n)$,
tend to the moments of the corresponding Gaussian pseudomatrix $\omega$, namely
$$
\lim_{n\rightarrow \infty}\tau_1(n)\left(\left(Y(n)\right)^{m}\right)=\Psi_1((\omega)^{m})
$$
where 
$$
\omega=\omega_{2,1}+\omega_{2,2}
$$
since $\omega_{1,2}=\omega_{1,1}=0$ and that is why they do not appear in the above formula 
(each $\omega_{p,q}$ is associated with the scalar $b_{p,q}=d_{p}v_{p,q}$ and we have $d_1=0$). 
In the random matrix contex, this means that the sequence $(T_{1,1}(n))$ is 
evanescent and $(T_{1,2}(n))$ is unbalanced. Moreover, 
$$
b_{2,1}=d_2v_{2,1}:=\beta_1\;\;{\rm and}\;\;b_{2,2}=d_2v_{2,2}:=\beta_2
$$ 
since $d_2=1$ and $v_{2,1}=v_{1,2}$. This proves the assertion in the case when 
$\alpha_1=\alpha_2=0$ (this includes Kesten laws). 

Before we prove the assertion for the general case, let us observe that the block refinement of the above asymptotics can be written in the form
$$
\lim_{n\rightarrow \infty}\tau_1(n)(T_{p_1,q_1}T_{p_2,q_2}\ldots T_{p_m,q_m})=
\Psi_1({\omega}_{p_1,q_1}\omega_{p_2,q_2}\ldots \omega_{p_m,q_m})
$$
provided we denote by $T_{2,1}$ rather than by $T_{1,2}$ the off-diagonal symmetric block.
Namely, by [11, Theorem 5.1], the mixed moments of symmetric blocks $T_{p,q}$ under partial traces converge 
to the corresponding mixed moments of {\it symmetrized} Gaussian operators $\widehat{\omega}_{p,q}$, where 
$\widehat{\omega}_{1,1}=\omega_{1,1}$  and $\widehat{\omega}_{2,2}=\omega_{2,2}$ and, more importantly, 
$$
\widehat{\omega}_{1,2}=\omega_{1,2}+\omega_{2,1}.
$$
Since, in the case considered in this theorem, $\omega_{1,2}=d_1v_{1,2}=0$ and thus 
$\widehat{\omega}_{1,2}=\omega_{2,1}$, we can replace each $\widehat{\omega}_{p_i,q_i}$ by 
$\omega_{p_i,q_i}$, which leads to the above equation. Moreover, even more information about
these moments can be obtained. For that purpose, decompose 
${\mathbb C}^{n}=W_1\oplus W_2$, where $W_j$ is the linear span of basis vectors 
indexed by $i\in N_j$ and observe that
$$
T_{2,1}(W_1)\subseteq W_2, \;
T_{2,1}(W_2)\subseteq W_1\;\;{\rm and}\;\;T_{j,j}(W_j)\subseteq W_j
$$
for $j\in \{1,2\}$. Since $\tau_1(n)$ is the partial trace over basis vectors from $W_1$, 
the above mixed moments of symmetric blocks vanishes unless it takes the form
in which even powers of $T_{2,2}$ alternate with $T_{2,1}$, namely
$$
\tau_1(n)(T_{2,1}T_{2,2}^{m_1}T_{2,1}\ldots 
T_{2,1}T_{2,2}^{m_r}T_{2,1}),
$$
where $m_1, \ldots , m_r\in 2{\mathbb N}\cup \{0\}$ and $m_1+m_2+\ldots +m_r+2r=m$. 
Likewise, the corresponding mixed moments of 
matricially free Gaussian operators vanish unless they take the form
$$
\Psi_1({\omega}_{2,1}\omega_{2,2}^{m_1}\omega_{2,1}\ldots {\omega}_{2,1}
\omega_{2,2}^{m_r}\omega_{2,1})
$$
since $\omega_{2,1}$ acts non-trivially onto $\Omega_1$ giving $e_1$ and sends $e_1$ back
to $\Omega_1$, whereas $\omega_{2,2}$ kills both $\Omega_1$ and $e_1$, leaving ${\mathcal F}_{2}$ invariant.
An even more detailed inspection leads to the formula
$$
\lim_{n\rightarrow \infty}\tau_1(n)(S_{1,2}T_{2,2}^{m_1}S_{2,1}\ldots 
S_{1,2}T_{2,2}^{m_r}S_{2,1})=
\Psi_1(\wp_{2,1}^{*}\omega_{2,2}^{m_1}\wp_{2,1}\ldots \wp_{2,1}^{*}\omega_{2,2}^{m_1}\wp_{2,1})
$$
since $\wp_{2,1}\Omega_1=e_1$ and $\wp_{2,1}^{*}e_1=\Omega_1$. Note that
the last formula is not obvious since it is not true in general that $S_{1,2}\rightarrow \wp_{2,1}^{*}$ and 
$S_{2,1}\rightarrow \wp_{2,1}$ under the partial traces.
However, it is very convenient because it allows us to study the effect of inserting the diagonal deterministic matrix 
$$
B=\alpha_1I_1(n) +\alpha_2I_2(n)
$$
between the symmetric blocks, where the dependence of $B$ on $n$ is supressed. 
We will show that an insertion of $B$ somewhere on the LHS of the above formula corresponds 
to an insertion of the operator $\gamma$ at the corresponding place on the RHS. Namely, this local analysis gives:
\begin{enumerate}
\item
at the left or right end of the above moment, the matrix $B$ reduces to $\alpha_1I_{1}$ and thus it produces 
$\alpha_1$ since it acts onto $W_1$; the corresponding $\gamma$ can also be replaced by $\alpha_1$ since it acts onto $\Omega$,
\item
in products of type $BS_{2,1}$ and $BT_{2,2}$, the matrix $B$ reduces to $\alpha_2I_2$ and gives $\alpha_2$ since it acts onto $W_2$; the corresponding pairs $\gamma\wp_{2,1}$ and $\gamma\omega_{2,2}$ can be replaced by 
$\alpha_2\wp_{2,1}$ and $\alpha_{2}\omega_{2,2}$, respectively, 
since $\gamma$ acts here onto vectors from ${\mathcal F}_2$.
\end{enumerate}
Consequently, for all non-trivial mixed moments of $T_{p,q}$ and $B$, we can write
$$
\lim_{n\rightarrow \infty}\tau_1(n)(B^{n_0}YB^{n_1}Y\ldots YB^{n_k})
=
\Psi_1(\gamma^{n_0}\omega\gamma^{n_1}\omega\ldots \omega\gamma^{n_k})
$$
for any nonnegative integers $n_0,n_1, \ldots, n_k$ and any 
$\alpha_1$ and $\alpha_2$. This implies that
$$
\lim_{n\rightarrow \infty}\tau_1(n)\left((M(n))^{m}\right)=\Psi_1((\omega+\gamma)^{m}),
$$
which completes the proof of our theorem.
\hfill $\blacksquare$

\begin{Corollary}
If $\beta_1=v_{2,1}>0$ and $\beta_2=v_{2,2}=0$ and under the remaining assumptions as in Theorem 4.1, it holds that
$$
\lim_{n\rightarrow \infty}\tau_1(n)\left((M(n))^{m}\right)=
\Psi_1((\omega_1+\gamma_1)^{m})
$$
where $\omega_1, \gamma_1$ are given by Corollary 3.1.
\end{Corollary}
{\it Proof.}
The proof is similar to that of Theorem 4.1. The only difference 
is that blocks $T_{2,2}(u,n)$ disappear from the computations under the trace $\tau_1(n)$
and thus non-trivial mixed moments take the special form
$$
\tau_1(n)(B^{n_0}T_{2,1}B^{n_{1}}T_{2,1}\ldots T_{2,1}B^{n_m})
=
\tau_1(n)(B^{n_0}S_{1,2}B^{n_{1}}S_{2,1}\ldots S_{2,1}B^{n_m})
$$
where $m$ is even and $S_{1,2}$ alternates with $S_{2,1}$. They tend to
$$
\Psi_1(\gamma^{n_0}\omega_{1}\gamma^{n_1}\omega_{1}\ldots \omega_{1}\gamma^{n_m})
=
\Psi_1(\gamma^{n_0}\wp_{1}^{*}\gamma^{n_1}\wp_{1}\ldots \wp_{1}\gamma^{n_m})
$$
as $n\rightarrow \infty$, where $\wp_1^{*}$ alternates with $\wp_{1}$, since each
$B^{j}S_{1,2}B^{k}$ can be replaced by $\alpha_{1}^{j}\alpha_2^{k}S_{1,2}$ for any $j,k\in {\mathbb N}$
by the definition of $B$ and, similarly, each $\gamma^{j}\wp_{1}^{*}\gamma^{k}$ can be replaced by
$\alpha_{1}^{j}\alpha_{2}^{k}\wp_{1}^{*}$ be the definition of $\gamma$. It remains to observe that
in the situation when we have mixed moments of $\omega_{1}$ and $\gamma$ under $\Psi_1$, we 
remain within ${\mathcal H}_{1}\oplus {\mathbb C}\Omega_1$ and thus $\gamma$ can be repleced by $\gamma_1$, which 
completes the proof. \hfill $\blacksquare$

\begin{Corollary}
Under the assumptions of Theorem 4.1, it holds that
$$
\lim_{n\rightarrow \infty}\tau_2(n)\left((M(n))^{m}\right)=
\Psi_2((\omega_2+\gamma_2)^{m})
$$
where $\omega_2, \gamma_2$ are given by Corollary 3.2.
\end{Corollary}
{\it Proof.}
The proof is similar to that of Theorem 4.1. In this case, when we compute the moments of $M(n)$
under $\tau_2(n)$, the mixed moments of $T_{1,1}(n),T_{2,1}(n),T_{2,2}(n)$ and $B$ become zero as $n\rightarrow \infty$ if 
there is $T_{1,1}(n)$ or $T_{2,1}(n)$ among them. On the level of matrices, this can be explained as follows: the fact that $(T_{2,1}(n))$ is unbalanced and 
is forced to act onto 'many' (of order $n)$ basis vectors from $W_2$ giving 'few' (of order smaller than $n$) 
basis vectors from $W_1$ makes the moment containing $T_{2,1}(n)$ vanish in the limit $n\rightarrow \infty$ (in other words, 
zero asymptotic dimensions cannot be associated with inner blocks). Of course, the case of $T_{1,1}(n)$ is clear since 
it is evanescent. On the operatorial level, the effect of this is that the moments involving $\omega_{1}$ do not contribute to the limit moments since all operators act within 
${\mathcal M}_{2}$, where $\omega_1$ is trivial and thus these moments reduce to the moments of $\omega_{2}$ and 
$\gamma$ under $\Psi_2$. Moreover, it is not hard to see that in fact $\gamma$ can be replaced with
$\gamma_2$, which is the restriction of $\gamma$ to ${\mathcal M}_{2}$.
\hfill $\blacksquare$\\

\section{Free Meixner Ensemble}

Let us consider an ensemble of independent random matrices of type 
considered in Section 4 and study their limit joint 
distributions under the state $\Psi_1$ as $n\rightarrow \infty$. 
The situation parallels that for the case of independent Gaussian random matrices and their asymptotic freeness [17]. 
As in Section 4, we will rely on the result derived in [12].

\begin{Definition}
{\rm By the {\it Free Meixner Ensemble} we will understand the family 
of independent $n\times n$ Hermitian Gaussian random matrices 
$\{M(u,n):n\in {\mathbb N},\,u\in\mathpzc{U}\}$, where matrices 
$$
M(u,n)=Y(u,n)+\alpha_1(u)I_1(n)+\alpha_2(u)I_2(n)
$$
satisfy the assumptions of Theorem 4.1 or Corollary 4.1 for any $u\in \mathpzc{U}$, where $\mathpzc{U}$ is an index set, 
with the constants $\alpha_1(u)$, $\alpha_2(u)$ as well as variances 
$\beta_1(u)=v_{2,1}(u)$, $\beta_2(u)=v_{2,2}(u)$ depending on $u\in \mathpzc{U}$.
In particular, we assume that all matrices are decomposed into blocks in the same 
fashion for any fixed $n$ and that their asymptotic dimensions are $d_1=0$ and $d_2=1$ for all 
$u$. }
\end{Definition}

We already know from Theorem 4.1 that the asymptotic distribution of $M(u,n)$ under the partial trace 
$\tau_1(n)$ is the free Meixner distribution associated with 
$$
(\alpha_1(u), \alpha_2(u), \beta_1(u), \beta_2(u)),
$$
but we would like to find an asymptotic relation between independent random matrices from this ensemble. 
This relation is expected to be of asymptotic freeness type. In fact, we will demonstrate that
the Free Meixner Ensemble is asymptotically conditionally free. As in Section 4, we exclude the case when $\beta_1(u)=0$ for some $u$ since 
in this case the corresponding matrix realization would be purely deterministic, but one can easily extend 
all results to include this case.

We also know from [12] that the Hermitian Symmetric Gaussian Block Ensemble 
$$
\{T_{p,q}(u,n): u\in \mathpzc{U}, n\in {\mathbb N}\}
$$
is asymptotically symmetrically matricially free, where
symmetric matricial freeness is a symmetrized
version of matricial freeness. More precisely, its asymptotics is determined by
operators of type $\widehat{\omega}_{p,q}(u)$ which are limit realizations of the corresponding 
symmetric blocks $T_{p,q}(u)$. 
We shall use the results of [12], where we also studied the family of their sums
$$
Y(u,n)=\sum_{p\leq q}T_{p,q}(u,n),
$$
in order to find the limit distributions of 
the Free Meixner Ensemble.

We used the mutlivariate matricially free Fock space of tracial type. 
The definition of ${\mathcal M}$ remains the same as in Section 3, but instead of one-dimensional Hilbert spaces, 
we take direct sums 
$$
{\mathcal H}_j=\bigoplus_{u\in \mathpzc{U}}{\mathcal H}_{j}(u)
$$
where ${\mathcal H}_{j}(u)={\mathbb C}e_{j}(u)$ for any $j\in \{1,2\}$ and $u\in \mathpzc{U}$, where
$\{e_j(u):j\in \{1,2\},\,u\in \mathpzc{U}\}$ is an orthonormal set.
Let 
$$
\mathpzc{B}=\{\Omega_1, \Omega_2, e_{2}(u_1, \ldots , u_n),\, 
e_{2}(u_1, \ldots , u_{n-1})\otimes e_{1}(u_n): u_1, \ldots , u_n\in \mathpzc{U}, \,n\in \mathbb{N}\}
$$
be the orthonormal basis of ${\mathcal M}$, where we use a shorthand notation
$$
e_{2}(u_1, \ldots , u_n)=e_{2}(u_1)\otimes \ldots \otimes e_{2}(u_n).
$$
Then we define the family 
of creation operators $\wp_1(u), \wp_2(u)$ by the following rules:
\begin{eqnarray*}
\wp_{1}(u)\Omega_1&=&\sqrt{\beta_{1}(u)}\,e_{1}(u)\\
\wp_{2}(u)\Omega_2&=&\sqrt{\beta_{2}(u)}\,e_{2}(u)\\
\wp_{2}(u)e_{2}(u_1,\ldots, u_n)
&=&
\sqrt{\beta_{2}(u)}\,e_{2}(u,u_1,\ldots, u_{n})\\
\wp_{2}(u)e_{2}(u_1,\ldots , u_{n-1})\otimes e_{1}(u_n)
&=&
\sqrt{\beta_{2}(u)}\,e_2(u,u_1,\ldots, u_{n-1})\otimes e_{1}(u_n)
\end{eqnarray*}
and we assume that $\wp_{1}(u), \wp_2(u)$ send the remaining basis vectors to zero. 
By $\wp_{1}^{*}(u)$ and $\wp_{2}^{*}(u)$ we denote their adjoints, respectively, 
and sums of the form
$$
\omega_{j}(u)=\wp_{j}(u)+\wp_{j}^{*}(u)
$$
are the corresponding Gaussian operators. We have shown in [12] that operators of this type give the limit realization of the mixed moments of symmetric blocks of independent Hermitian Gaussian random matrices with block-identical variances (Gaussian Symmetric Block Ensemble). 
In other words, we showed that we have convergence of mixed moments
$$
\lim_{n\rightarrow \infty}\tau_q(n)(T_{p_1,q_1}(u_1,n)\ldots T_{p_m,q_m}(u_m,n))
=
\Psi_q(\widehat{\omega}_{p,q}(u_1)\ldots \widehat{\omega}_{p_m,q_m}(u_m)).
$$
where $\widehat{\omega}_{p,q}(u)$ is the same symmetrization as in the case of $\widehat{\omega}_{p,q}$ 
in Section 3.

Let us give a definition of conditional freeness which is very similar to that of freeness and that
will be helpful for us. The family of unital subalgebras $\{{\mathcal A}(u): u\in \mathpzc{U}\}$ of a unital algebra ${\mathcal A}$ 
is {\it conditionally free} with respect to the pair of states $(\varphi, \psi)$ on ${\mathcal A}$ if 
$$
\varphi(a_1a_2\ldots a_m)=0
$$
whenever $a_i\in {\mathcal A}(u_i)\cap {\rm Ker}\psi$ for any $1\leq i \leq n-1$ and 
$a_n\in {\mathcal A}(u_n)\cap {\rm Ker}\varphi$, where $u_1\neq u_2\neq \ldots \neq u_n$.
This definition is equivalent to other definitions and immediately shows that there is 
a relation between different levels of Hilbert spaces in their free product and the corresponding states 
assigned to these levels. Consequently, there is a relation 
with the depths of the blocks of noncrossing partitions which contribute to the moments of conditionally free 
random variables. In more generality, we obtain freeness with infinitely many states [7].

\begin{Theorem}
Let $\tau_j(n)$ be the partial trace over the set of basis vectors indexed by $N_j$, where $j\in \{1,2\}$. 
The family of matrices
$$
\{M(u,n):u\in \mathpzc{U},n\in {\mathbb N}\}
$$ 
is asymptotically conditionally free with respect to the pair of partial traces
$(\tau_1(n), \tau_2(n))$ as $n\rightarrow \infty$.
\end{Theorem}
{\it Proof.}
In particular, if we consider the $2\times 2$ block random matrices 
with asymptotic dimensions $d_1=0$ and $d_2=1$, the sequence $(T_{1,1}(n,u))$ is evanescent and $(T_{1,2}(n,u))$
is unbalanced and thus the corresponding arrays of symmetrized Gaussian operators reduce to
arrays containing only $\omega_2(u)=\omega_{2,2}(u)$ and $\omega_1(u)=\omega_{2,1}(u)$ simply because $\omega_{1,1}(u)=0$ and $\omega_{1,2}(u)=0$.
Thus, in view of the above, we have
$$
\lim_{n\rightarrow \infty}\tau_1(n)(Y(u_1,n)\ldots Y(u_m,n))
=
\Psi_1(\omega(u_1)\ldots \omega(u_m)),
$$
where
$$
\omega(u)=\omega_1(u)+\omega_2(u)
$$
for any $u\in \mathpzc{U}$. 
As in the proofs of Theorems 3.1 and 4.1, this can be generalized to the moments of
matrices $M(u,n)$ from the Free Meixner Ensemble since all computations presented there
are based on the relations between matricial indices of the considered blocks and of the considered
operators and they depend on $u$ only in the sense that the blocks associated with symmetric blocks
and with the corresponding operators labelled by $u$ give rise to parameters $\alpha_j(u), \beta_j(u)$
labelled by $u$. Thus, we have
$$
\lim_{n\rightarrow \infty}\tau_1(n)(M(u_1,n)\ldots M(u_m,n))
=
\Psi_1(y(u_1)\ldots y(u_m)),
$$
where 
$$
y(u)=\omega(u)+\gamma(u)
$$
and 
$$
\gamma(u)=(\alpha_{2}(u)-\alpha_1(u))(\beta_{1}^{-1}(u)\wp_{1}(u)\wp_{1}^{*}(u)+
\beta_{2}^{-1}(u)\wp_{2}(u)\wp_{2}^{*}(u))+\alpha_{1}(u),
$$
for any $u\in \mathpzc{U}$, where $\alpha_q(u)\in {\mathbb R}$ and $\beta_q(u)>0$ for $q\in \{1,2\}$.
In a similar way one shows that 
$$
\lim_{n\rightarrow \infty}\tau_2(n)(M(u_1,n)\ldots M(u_m,n))
=
\Psi_2(y(u_1)\ldots y(u_m))
$$
where Corollaries 3.2 and 4.2 are used. Therefore, in order to prove our assertion, we need to show that the family
$$
\{y(u):\,u\in \mathpzc{U}\}
$$
is conditionally free with respect
to the pair of states $(\Psi_1, \Psi_2)$, where $\Psi_q$ is the vector 
state associated with $\Omega_q$. We will prove a slightly more general result, namely that
the family of unital *-algebras $\{{\mathcal A}(u):u\in \mathpzc{U}\}$, each 
generated by $\wp_{2,1}(u)$ and $\wp_{2,2}(u)$ for fixed $u$, respectively, 
is conditionally free with respect to $(\Psi_1, \Psi_2)$. We need to show that
$$
\Psi_1(a_1a_2\ldots a_n)=0
$$
for any $a_i\in {\mathcal A}(u_i)\cap {\rm Ker}\Psi_2$, where $1\leq i \leq n-1$ and
$a_n\in {\mathcal A}_{u_n}\cap {\rm Ker}\Psi_1$. 

We claim that the variable $a_n$ is a polynomial in noncommuting variables 
$$
\wp_{1}(u_n),\wp_{1}^{*}(u_n), \wp_2(u_n), \wp_2^{*}(u_n)
$$ 
which can be written as a linear combination of $P^{\perp}=1-P$ and of monomials 
$$
\wp_{2}^{m_2}(u_n)\wp_{1}^{m_1}(u_n)(\wp_{1}^{*}(u_n))^{k_2}(\wp_{2}^{*}(u_n))^{k_2}
$$
where $m_2,k_2\in {\mathbb N}\cup \{0\}$, $k_1,k_2\in \{0,1\}$ are such that $m_1+m_2+k_1+k_2>0$.
In order to reduce all monomials from ${\rm Ker}\Psi_1$ to this form, first observe that 
${\mathcal M}_{1}$ is invariant under the action of $\wp_{1},\wp_2$ and their adjoints. 
Therefore, it suffices to consider all operators as their restrictions to ${\mathcal M}_{1}$. Then 
we have the relations
$$
\wp_{1}^{*}(u)\wp_1(u)=\beta_{1}(u)P, \;\;
\wp_{2}^{*}(u)\wp_2(u)=\beta_{2}(u)P^{\perp}
\;\;{\rm and}\;\;\wp_{q}^{*}(u)\wp_q(u')=0
$$
as well as 
$$
\wp_{1}(u)\wp_{2}(u')=0,
P\wp_{1}(u)=0, P^{\perp}\wp_{1}(u)=\wp_{1}(u), \wp_{1}(u)P=\wp_{1}(u), \wp_{1}(u)P^{\perp}=0
$$
for any $q\in \{1,2\}$ and $u\neq u'$, as well as their adjoints. 
Clearly, $P^{\perp}\in {\rm Ker}\Psi_1$.
Therefore, we can pull all starred operators to the right of the unstarred ones 
in $a_n$ and our claim is proved. This implies that $a_n$ maps $\Omega_1$ into 
${\mathcal M}_{1}\ominus {\mathbb C}\Omega_1$. 

Now, any vector from the image 
$a_n(\Omega_1)$ is a linear combination of vectors which begin with $e_2(u_n)$. 
Therefore, the action of $a_{n-1}$ onto these vectors is the same as its
action onto $\Omega_2$. Therefore, if we take $a_{n-1}\in {\rm Ker}\Psi_2$ 
and we apply a similar reasoning as above, we can write it as a linear combination of
monomials
$$
\wp_{2}^{m_2}(u_{n-1})(\wp_{2}^{*}(u_{n-1}))^{k_2}
$$
where $m_2+k_2>0$ since ${\mathcal A}(u_{n-1})$ leaves ${\mathcal M}_{1}\ominus {\mathbb C}\Omega_1$
invariant (recall that $u_{n-1}\neq u_n$) and the action of $\wp_{1}^{*}(u_{n-1})$ is trivial
on this space.
Moreover, the constant term vanishes since $a_{n-1}\in {\rm Ker}\Psi_2$ and thus
$a_{n-1}a_{n}(\Omega_1)$ is a linear combination of vectors which begin with $e_{2}(u_{n-1})$.
Continuing in this fashion, we obtain $a_1a_2\ldots a_n(\Omega_1)\perp \Omega_1$, which completes the proof.
\hfill $\blacksquare$

\begin{Remark}
{\rm It can be easily seen that the family of algebras $\{{\mathcal A}(u):u\in \mathpzc{U}\}$ is, in general, not free
with respect to $\Psi_1$. For instance, in the simple case when the Jacobi parameters are $(0,0, \beta_1, \beta_2)$
for all $u\in \mathpzc{U}$, where $0\neq \beta_{1}\neq \beta_{2}\neq 0$, then we can take two polynomials, say 
$w_1=y(s)\in {\mathcal A}(s)$ and $w_2=(y(u))^{2}-\beta_{1}\in {\mathcal A}(u)$, where $u\neq s$, which are 
in ${\rm Ker}\Psi_1$, but 
$$
\Psi_1(w_1w_2w_1)=\beta_{1}(\beta_{2}-\beta_{1})\neq 0
$$
since 
$$
\wp_{1}^{*}(u)\wp_{1}(u)\Omega_1=\beta_{1}\Omega_1\;\;{\rm and}\;\; \wp_{2}^{*}(u)\wp_{2}(u)e_1=\beta_{2}e_1.
$$ 
Of course, if we replace $w_2$ by $w_{3}=(y(u))^{2}-\beta_{2}\in {\rm Ker}\Psi_2$, we get zero
in the above equation, which in in agreement with the conditional freeness of $w_1,w_3$ with 
respect to $(\Psi_1, \Psi_2)$ stated in Theorem 5.1.}
\end{Remark}

\end{document}